\documentclass[12pt]{amsart}

\usepackage{amsmath,amssymb,amscd,amsthm,fullpage}
\usepackage[mathscr]{eucal}

\newtheorem{theorem}{Theorem}[section]
\newtheorem{proposition}[theorem]{Proposition}
\newtheorem{corollary}[theorem]{Corollary}
\newtheorem{lemma}[theorem]{Lemma}

\theoremstyle{definition}

\newtheorem*{definition*}{Definition}

\newcommand{\Z}{\mathbb{Z}}

\newcommand{\Q}{\mathbb{Q}}
\newcommand{\eps}{\varepsilon}
\newcommand{\p}{\partial}
\newcommand{\half}{\tfrac{1}{2}}
\newcommand{\CP}{\mathbb{CP}}
\renewcommand{\*}{\cdot}
\newcommand{\E}{\mathsf{E}}
\newcommand{\PP}{\mathsf{P}}

\DeclareMathOperator{\Res}{Res}
\DeclareMathOperator{\res}{res}

\renewcommand{\]}{{]\!]}}

\newcommand{\Om}{\Omega}

\renewcommand{\o}{\otimes}
\renewcommand{\)}{{)\!)}}

\DeclareMathOperator{\ad}{ad}
\renewcommand{\AA}{\mathcal{A}}
\newcommand{\tAA}{\widetilde{\mathcal{A}}}
\newcommand{\BB}{\mathcal{B}}
\renewcommand{\d}{\delta}
\newcommand{\dbar}{\bar\d}
\newcommand{\Lbar}{\bar{L}}
\newcommand{\Wbar}{\overline{W}}
\newcommand{\abar}{\bar{a}}
\newcommand{\pbar}{\bar{p}}
\newcommand{\qbar}{\bar{q}}
\newcommand{\vbar}{\bar{v}}
\newcommand{\xbar}{\bar{x}}
\newcommand{\zbar}{\bar{z}}
\newcommand{\Hbar}{\bar{H}}
\renewcommand{\t}{\nu}
\newcommand{\tint}{{\textstyle\int}}

\newcommand{\RR}{\mathcal{R}}

\begin{document}

\title{The equivariant Toda lattice, II}

\author{Ezra Getzler}

\address{Department of Mathematics, Northwestern University, Evanston,
  IL, 60208, USA}

\maketitle

This paper is a continuation of \cite{eqtoda}, in which we introduced a
reduction of the Toda lattice hierarchy (in the limit of infinitesimal
lattice spacing), called the equivariant Toda lattice, by imposing the
constraints
\begin{align} \label{Equivariant}
(\d_1 - \dbar_1)L &= \t\p L , & (\d_1 - \dbar_1)\Lbar &= \t\p\Lbar
\end{align}
on the Lax operators $L$ and $\Lbar$. This reduction is a deformation of
the the Toda chain, which is the reduction corresponding to the constraint
$L=\Lbar$.

Seeking an integrable system which would describe the Gromov-Witten
invariants of $\CP^1$, Eguchi and Yang \cite{EY} studied the Toda
chain. They conjectured the existence in the limit of infinitesimal
lattice spacing of an additional hierarchy of commuting flows: these
flows were constructed independently by the author \cite{toda} and
Zhang \cite{Zhang} using homological perturbation theory, and shown to
be bihamiltonian. Recently, Carlet, Dubrovin and Zhang \cite{CDZ} have
shown that these additional flows may be described by Lax equations
involving the logarithm of the Lax operator $L$, as conjectured by
Eguchi and Yang.

In this paper, borrowing the ideas of Carlet et al., we prove that the
equivariant Toda lattice has a Hamiltonian structure which is a deformation
of the first Hamiltonian structure of the Toda chain. (We were however
unable to find a bihamiltonian structure.)

\subsection*{Acknowledgements}

This paper was inspired by discussions with B. Dubrovin and Y. Zhang at the
meeting on Frobenius manifolds at the Max-Planck-Institut f\"ur Mathematik,
Bonn, and the author thanks C. Hertling, Yu.\ Manin and M. Marcolli for the
invitation to participate. The author is partially supported by the NSF
under grant DMS-0072508.

\section{Notation}

In this section, we recall some of the terminology of \cite{eqtoda}. All of
the commutative algebras which we consider carry an involution
$p\mapsto\pbar$. By a \textbf{differential algebra}, we mean a commutative
algebra with derivation $\p$ such that
$$
\p\pbar = \overline{\p p} .
$$
If $\AA$ is a differential algebra and $S$ is a set, the free
differential algebra $\AA\{S\}$ generated by $S$ is the polynomial
algebra
$$
\AA[\p^nx,\p^n\xbar\mid x\in S, n\ge0] ,
$$
with differential $\p(\p^nx)=\p^{n+1}x$. An \textbf{evolutionary
derivation} $\delta$ of a differential algebra $\AA$ is a derivation such
that $[\p,\delta] = 0$.

Let $\AA$ be a differential algebra over $\Q_\eps$, and let $q\in\AA$ be a
regular element (that is, having no zero-divisors) such that $\qbar=q$. The
localization $q^{-1}\AA$ of $\AA$ is a differential algebra, with
differential $\p(q^{-1})=-q^{-2}\p q$.  Let $\Phi_\pm(\AA,q)$ be the
associative algebras of \textbf{difference operators}
\begin{align*}
\Phi_+(\AA,q) &= \biggl\{ \sum_{k=-\infty}^\infty p_k \, \Lambda^k \biggm|
\text{$p_k\in q^{-1}\AA$, $p_k=0$ for $k\ll0$} \biggr\} , \\
\Phi_-(\AA,q) &= \biggl\{ \sum_{k=-\infty}^\infty p_k \, \Lambda^k \biggm|
\text{$p_k\in\AA$, $p_k=0$ for $k\gg0$} \biggr\} ,
\end{align*}
with product
$$
\sum_i a_i \Lambda^i \* \sum_j b_j \Lambda^j = \sum_k \biggl( \sum_{i+j=k}
\bigl(\E^{-j/2}a_i\bigr)\bigl(\E^{i/2}b_j\bigr) \biggr) \Lambda^k .
$$
Note that $\Phi_-(\AA,q)$ is in fact independent of $q$.

Let $A\mapsto A_\pm$ be the projections from on $\Phi_\pm(\AA,q)$ defined
by the formulas
\begin{align*}
\biggl( \sum_{k=-\infty}^\infty p_k\Lambda^k \biggr)_+ &= \sum_{k=0}^\infty
p_k\Lambda^k , &
\biggl( \sum_{k=-\infty}^\infty p_k\Lambda^k \biggr)_- &=
\sum_{k=-\infty}^{-1} p_k\Lambda^k .
\end{align*}
We see that $A=A_-+A_+$. Define the residue $\res:\Phi_\pm(\AA,q)\to\AA$ by
the formula
$$
\res\biggl( \sum_{k=-\infty}^\infty p_k\Lambda^k \biggr) = p_0 .
$$
For $k\in\Z$, let $[k]$ be the isomorphism of $\AA$
$$
[k] = \frac{\E^{k/2}-\E^{-k/2}}{\E^{1/2}-\E^{-1/2}} =
\sum_{j=1}^k \E^{(k+1)/2-j} = k + O(\eps^2) .
$$
Define $q^{[k]}$ by the recursion
$$
q^{[k+1]} = \E^kq \* \E^{-1/2}q^{[k]} ,
$$
with initial condition $q^{[0]}=1$. The involution
$$
A = \sum_{k=-\infty}^\infty p_k \Lambda^k \mapsto \bar{A} =
\sum_{k=1}^\infty \pbar_k \, q^{[k]} \Lambda^{-k} + \pbar_0 +
\sum_{k=1}^\infty \pbar_{-k} \, q^{-[k]} \Lambda^k ,
$$
defines an anti-isomorphism between the algebras $\Phi_+(\AA,q)$ and
$\Phi_-(\AA,q)$.

\section{The dressing operator of the Toda lattice}

Let $\BB$ be the free differential algebra $\Q_\eps\{q,w_k\mid
k>0\}/(q-\qbar)$, and let $W$ be the universal \textbf{dressing operator}
of the Toda lattice
$$
W = 1 + \sum_{k=1}^\infty w_k \Lambda^{-k} \in \Phi_-(\BB,q) .
$$
The coefficients $w_k^*\in\BB$ of $W^{-1}$,
$$
W^{-1} = 1 + \sum_{k=1}^\infty w_k^* \Lambda^{-k} ,
$$
are characterized by the recursion
$$
w_k^* = - w_k - \sum_{j=1}^{k-1}
\bigl(\E^{(k-j)/2}w_j\bigr)\bigl(\E^{-j/2}w_{k-j}^*\bigr) ,
$$
obtained by extracting the coefficient of $\Lambda^{-k}$ in the equation
$WW^{-1}=I$.

The Lax operator of the Toda lattice is the difference operator
$$
L = W\Lambda W^{-1} = \Lambda + \sum_{k=1}^\infty a_k\Lambda^{-k+1} \in
\Phi_-(\BB,q) .
$$
Since $a_k+\eps\nabla w_k$ lies in the differential ideal
$(w_1,\dots,w_{k-1})$ for all $k>0$, we see that the sequence of elements
$a_k$ of $\BB$ defines an embedding of differential algebras
$$
\AA=\Q_\eps\{q,a_k\mid k>0\}/(q-\qbar) \hookrightarrow \BB .
$$

The conjugate Lax operator $\Lbar$ is
$$
\Lbar = \Wbar^{-1}(q\Lambda^{-1})\Wbar = q\Lambda^{-1} + \sum_{k=1}^\infty
\abar_k q^{-[k-1]} \Lambda^{k-1} .
$$
Let $B_n=\eps^{-1}L^n_+$ and $C_n=-\eps^{-1}\Lbar^n_-$.

We define evolutionary derivations $(\d_n,\dbar_n\mid n>0)$ of $\BB$ by the
formulas
\begin{equation} \label{flow}
\eps\d_nW + L^n_-W = \eps\dbar_nW + \Lbar^n_-W = 0 .
\end{equation}
These derivations are called the flows of the Toda lattice. The action of
the derivations $\d_n$ and $\dbar_n$ on $\BB$ restricts to an action on
$\AA$ such that the derivatives of the Lax operator $L$ are given by the
Lax equations $\d_nL=[B_n,L]$ and $\dbar_nL=-[C_n,L]$. These flows on $\BB$
commute, by the Zakharov-Shabat equations
\begin{align*}
\d_mB_n - \d_nB_m &= [B_m,B_n] , &
\d_mC_n - \dbar_nB_m &= [B_m,C_n] , &
\dbar_mC_n - \dbar_nC_m &= [C_m,C_n] ,
\end{align*}
and $\d_n$ is indeed the conjugate derivation to $\dbar_n$. 

Let $\log(L)=W\log(\Lambda)W^{-1}$, where $\log(\Lambda)$ is a formal
symbol for the operator $\eps\p$: namely, we have the commutation relation
$$
[\log(\Lambda),f] = \eps \p f , \quad f\in\Phi_-(\AA,q) .
$$
Define $\ell$ to be the difference operator
\begin{align*}
\ell &= \log(\Lambda) - \log(L) = \eps(\p W)W^{-1} \\
&= \eps \biggl( \p w_k + \sum_{j=1}^{k-1}
\bigl(\E^{(k-j)/2}\p w_j\bigr)\bigl(\E^{-j/2}w_{k-j}^*\bigr) \biggr) .
\end{align*}
The following is a result of Carlet, Dubrovin and Zhang \cite{CDZ}.
(They work in the context of the Toda chain, so they assume that $a_k=0$,
$k>2$.)
\begin{proposition} \label{ell}
The difference operator $\ell$ is an element of $\Phi_-(\AA,q)$.
\end{proposition}
\begin{proof}
Write
$$
\ell = \sum_{k=1}^\infty b_k \Lambda^{-k} \in \Phi_-(\BB,q) .
$$
We show that $b_k\in\AA$ for all $k>0$, by induction on $k$.

Define elements $p_k(n)$ of $\AA$ by the formula
$$
L^n = \sum_{k=-\infty}^\infty p_k(n) \Lambda^k .
$$
We have
$$
\eps \p L = \eps \p (W\Lambda W^{-1})
= \eps (\p W)\Lambda W^{-1} - \eps W\Lambda W^{-1}(\p W)W^{-1} = [\ell,L] ,
$$
hence for each $n>0$, $\eps\p L^n=[\ell,L^n]$. Applying
the linear map $\res:\Phi_-(\AA,q)\to\AA$, we obtain the equation
\begin{equation} \label{bn}
\nabla \Biggl( [n]b_n + \sum_{k=1}^{n-1} [k] ( b_kp_k(n) ) + \PP
p_0(n) \Biggr) = 0 .
\end{equation}

Denote by $\alpha:\AA\to\Q_\eps$ the homomorphism which sends the
generators $\{q,a_k,\abar_k\}$ of $\AA$ to $0$. Since $\alpha\*\p=0$, we
see that $\alpha(\p W)=0$, and hence $\alpha(\ell)=0$. Thus, the constant
of integration in \eqref{bn} vanishes, and we obtain the recursive formula
\begin{equation} \label{Bn}
b_n = - \frac{1}{[n]} \Biggl( \sum_{k=1}^{n-1} [k] ( b_kp_k(n) ) +
\PP p_0(n) \Biggr)
\end{equation}
for the coefficients $b_k$, showing that they are elements of $\AA$.
\end{proof}

\section{Fractional powers of the Lax operator}

In this section, we study the fractional powers of the Lax operator $L$;
this may be compared with the parallel construction for the KP hierarchy
due to Khesin and Zakharevich \cite{KZ}. The study of these fractional
powers is closely related to the operator $\ell$ introduced in the last
section.

Let $s$ be a complex number. The fractional power $L^s$ of the Lax operator
$L$ is defined by means of the dressing operator:
\begin{equation} \label{ls}
L^s = W\Lambda^sW^{-1} = \Lambda^s + \sum_{k=1}^\infty a_k(s)
\Lambda^{s-k} \in \Phi_-(\BB,q) .
\end{equation}
The coefficient $a_k(s)$ is given by the explicit formula
$$
a_k(s) = \E^{-s/2}w_k + \sum_{j=1}^{k-1}
\bigl( \E^{(k-j-s)/2}w_j \bigr) \bigl( \E^{(s-j)/2}w_{k-j}^* \bigr) +
\E^{s/2}w_k^* .
$$
In particular, $a_k(0)=0$ and $a_k(1)=a_k$. Differentiating the definition
\eqref{ls} of $L^s$ with respect to $s$ and setting $s=0$, we obtain the
formula
\begin{equation} \label{Ls0}
\frac{dL^s}{ds} \Big|_{s=0} = - \ell ,
\end{equation}
showing that $a'_k(0)=-b_k$. The following proposition is proved by
extending this differential equation to all values $s$.
\begin{proposition} \label{aks}
The coefficient $a_{k,i}(s)$ in the expansion
$$
a_k(s) = \sum_{i=0}^\infty \eps^i a_{k,i}(s)
$$
is a polynomial in $s$ of degree $i+1$ with coefficients in the
differential algebra
$$
\Q\{q,a_k\mid k>0\}/(q-\qbar) .
$$
\end{proposition}
\begin{proof}
By its definition, the fractional power $L^s$ satisfies the differential
equation
$$
\frac{dL^s}{ds} = - \half \bigl( L^s\,\ell + \ell\,L^s \bigr) .
$$
Taking the coefficient of $\Lambda^{s-k}$ on both sides, we obtain the
differential equation
$$
\frac{da_k(s)}{ds} = - \frac12 \sum_{j=1}^{k-1} \bigl( \E^{(s-j)/2}b_{k-j} \,
\E^{(k-j)/2}a_j(s) + \E^{(j-s)/2}b_{k-j} \, \E^{(j-k)/2}a_j(s) \bigr) ,
$$
where we interpret $a_0(s)$ as $1$. By an application of Proposition
\ref{ell}, the result follows.
\end{proof}

\section{Perturbation theory}

Let $\Om(\AA)$ be the vector space of K\"ahler differentials of the
commutative $\Q_\eps$-algebra $\AA$; this is a free module over $\AA$ with
basis $\{dq,da_k,d\abar_k\mid k>0\}$. The differential $d:\AA\to\Om(\AA)$
extends to a morphism
$$
d:\Phi_-(\AA,q)\to\Phi_-(\AA,q)\o_\AA\Om(\AA) .
$$
The goal of this section is the calculation of the differentials $dL^s$ and
$d\ell$ in terms of the fundamental differential
$$
dL = \sum_{k=1}^\infty da_k \, \Lambda^{-k+1} .
$$

A basic formula of perturbation theory (Kumar \cite{Kumar}) says that
for $f(z)$ an analytic function of $z$,
$$
df(L) = \sum_{k=0}^\infty \frac{(-1)^k}{(k+1)!} \ad(L)^k(f^{(k+1)}(L)dL) .
$$
For $f(z)=z^s$, this becomes
\begin{equation} \label{dLs}
dL^s = \sum_{k=0}^\infty (-1)^k \tbinom{s}{k+1} \ad(L)^k(L^{s-k-1}dL) .
\end{equation}
We will now prove this formula directly.

For $s$ a natural number $n$, the right-hand side of \eqref{dLs} is a
finite sum, and the formula is then easily proved by induction on $n$: we
have
\begin{align*}
d(L^{n+1}) &= dL^n\*L + L^n\*dL = \sum_{k=0}^{n-1} (-1)^k \tbinom{n}{k+1}
\ad(L)^k(L^{n-k-1}dL)\*L + L^n \* dL \\
&= \sum_{k=0}^{n-1} (-1)^k \tbinom{n}{k+1}
\bigl( \ad(L)^k(L^{n-k}dL) - \ad(L)^{k+1}(L^{n-k-1}dL) \bigr) + L^n \* dL
\\
&= \sum_{k=0}^n (-1)^k \bigl( \tbinom{n}{k} + \tbinom{n}{k+1} \bigr)
\ad(L)^k(L^{n-k}dL)
= \sum_{k=0}^n (-1)^k \tbinom{n+1}{k+1} \ad(L)^k(L^{n-k}dL) .
\end{align*}

By analytic continuation, \eqref{dLs} holds for all values of $s$. Indeed,
the right-hand side is convergent in the $\eps$-adic topology, since the
operation $\ad(L)$ may be split into two terms: $\ad(\Lambda+a_1)=O(\eps)$,
and
$$
\sum_{k=2}^\infty \ad( a_k \Lambda^{-k+1} ) = O(\Lambda^{-1}) .
$$
It only remains to observe that by Theorem \ref{aks}, the coefficient of
$\eps^i$ in $da_{k,i}(s)$ is polynomial in $s$.

It is now straightforward to calculate $d\ell$: taking the derivative of
\eqref{dLs} with respect to $s$ and setting $s=0$, we see that
\begin{equation} \label{dlog}
d\ell = - \sum_{k=0}^\infty \frac{1}{k+1} \ad(L)^k(L^{-k-1}\,dL) .
\end{equation}

\section{The equivariant Toda lattice and $\ell$}

In this section, we denote the element $a_1\in\AA$ by $v$. Let $K$ be the
difference operator
$$
K = L_+ + \Lbar_- = \Lambda+v+q\Lambda^{-1} \in
\Phi_+(\AA,q)\cap\Phi_-(\AA,q) .
$$
In \cite{eqtoda}, we defined the equivariant Toda lattice by the
constraints
\begin{align} \label{equivariant}
\eps^{-1}[K,L] &= \t\p L , & \eps^{-1}[K,\Lbar] &= \t\p\Lbar ,
\end{align}
or equivalently, the constraints \eqref{Equivariant}. We showed that the
differential algebra associated to the equivariant Toda lattice is
isomorphic to
$$
\tAA = \Q_{\eps,\t}[z_k,\zbar_k\mid k>0]\{q,v,\vbar\}/(\t\p
q-\nabla(v-\vbar)) ,
$$
where $\Q_{\eps,\t}=\Q_\eps[\t]$, and the constants of motion $z_k$ are
the images of the elements
$$
p_{-1}(k) - qp_1(k) - \t\PP p_0(k) \in \AA
$$
under the natural quotient map from $\AA$ to $\tAA$.

Let $e$ be the derivation $\p_v+\p_{\vbar}$ of $\tAA$; then $e(K)=1$ and
\begin{equation} \label{puncture}
\biggl( L - \t + \sum_{k=1}^\infty z_k L^{-k} \biggr) e(L) = L .
\end{equation}
\begin{theorem} \label{main}
The constraint \eqref{equivariant} defining the equivariant Toda lattice is
equivalent to the identity
\begin{equation} \label{KL}
K = L + \t \ell - \sum_{k=1}^\infty \frac{z_k}{k} L^{-k} .
\end{equation}
The vanishing of the constants $z_k$ is equivalent to the constraint
\begin{equation} \label{equivariantW}
(\d_1-\dbar_1)W = \t\p W ,
\end{equation}
or equivalently, the equation $(\d_1-\dbar_1)=\t\p$ on the differential
algebra $\BB$.
\end{theorem}
\begin{proof}
Written in terms of $\ell$, \eqref{equivariant} becomes
$$
[K-\t\ell,L] = 0 .
$$
This is equivalent to the statement that
$$
K - \t\ell \in \Q_{\eps,\t}\(L^{-1}\) .
$$
It is not hard to see that
\begin{equation} \label{c}
K - \t\ell - L = \sum_{k=1}^\infty y_k L^{-k} \in
\Q_{\eps,\t}\[L^{-1}\] ;
\end{equation}
the constant term vanishes since, by definition, $\res(K)$ and $\res(L)$
equal $v$, while $\res(\ell)=0$.

It remains to identify the constants $y_k$. If $\d$ is an evolutionary
derivation of the differential algebra $\tAA$, \eqref{dlog} implies that
$$
\d\ell = - \sum_{k=0}^\infty \frac{1}{k+1} \ad(L)^k(L^{-k-1}\d L) .
$$
In particular, since $L$ commutes with $e(L)$, we see that $e(\ell) = -
L^{-1} e(L)$. Likewise, $e(L^{-k})=-kL^{-k-1}$. Applying the derivation $e$
to both sides of \eqref{c}, we see that
$$
1 = e(K) = e(L) \biggl( 1 - \t L^{-1} - \sum_{k=1}^\infty k y_k L^{-k-1}
\biggr) .
$$
It follows from \eqref{puncture} that $y_k=-z_k/k$.

We have
\begin{align*}
(K-L-\t\ell)W &= (L_++\Lbar_-)W - LW - \eps\t\p W \\
&= - L_-W + \Lbar_-W - \eps\t\p W = \eps(\d_1-\dbar_1-\t\p)W .
\end{align*}
Thus, the vanishing of the constants $z_k$ in \eqref{KL} is equivalent to
the constraint \eqref{equivariantW}.
\end{proof}

In \cite{eqtoda}, we conjectured that the equivariant Gromov-Witten
invariants of $\CP^1$ are described by the equivariant Toda lattice with
$z_k=0$, $k>0$. The results of this section show that this is true. By the
work of Okounkov and Pandharipande, the equivariant Gromov-Witten
invariants of $\CP^1$ are associated with a $\tau$-function of the Toda
lattice which satisfies $(\d_1-\dbar_1)\tau=\t\p\tau$. The dressing
operator $W$ corresponding to this $\tau$-function is given by the formula
$$
W = \tau^{-1} \exp\biggl( - \sum_{n=1}^\infty \frac{\d_n}{n\Lambda^n}
\biggr) \tau ;
$$
it follows that $W$ satisfies the equation $(\d_1-\dbar_1)W=\t\p W$.

The other part of the conjecture of \cite{eqtoda}, relating the equivariant
Gromov-Witten flows of $\CP^1$ to the flows of the equivariant Toda lattice,
is established by Okounkov and Pandharipande. Namely, if $\p_k=\p_{k,Q}$
and $\bar{\p}_k=\p_{k,Q}-\t\p_{k,P}$, then
\begin{align*}
\sum_{k=0}^\infty z^{k+1} \p_k &= \sum_{n=1}^\infty
\frac{z^n\delta_n}{(1+z\t)(2+z\t)\dots(n+z\t)} , \\
\sum_{k=0}^\infty z^{k+1} \bar{\p}_k &= \sum_{n=1}^\infty
\frac{z^n\dbar_n}{(1-z\t)(2-z\t)\dots(n-z\t)} .
\end{align*}
In particular, we see that the descendent flows $\p_{k,P}$ of the puncture
operator $P$ are given in the non-equivariant limit by the formula
\begin{equation} \label{pP}
\p_{k,P} = \lim_{\t\to0} \bigl( \tfrac{1}{(k+1)!} \t^{-1} ( \d_{k+1} -
\dbar_{k+1} ) - \tfrac{1}{k!} c_k (\d_k+\dbar_k) \bigr) ,
\end{equation}
where $c_k$ is the harmonic number $c_k = 1 + \tfrac{1}{2} + \dots +
\tfrac{1}{k}$.

\section{Hamiltonian structure}

In this section, we use Theorem \ref{main} to show that the equivariant
Toda lattice has a Hamiltonian structure.

Denote by $\RR$ the quotient $\tAA/\p\tAA$, and denote by $f\mapsto\tint
f\,dx$ the quotient map from $\tAA$ to $\RR$. The idea which this notation
is intended to represent is that an element of $\tAA$ is a density $f$,
whose associated functional $\tint f\,dx$ is obtained by integration with
respect to the space variable $x$.

Denote by $\Res$ the trace on $\Phi_-(\tAA,q)$ with values in $\RR$ given
by the formula
$$
\Res\biggl( \sum_{k=-\infty}^\infty f_k\Lambda^k \biggr) = \tint f_0 \, dx .
$$
Clearly, this map vanishes on total derivatives; to see that it vanishes on
commutators, we use the formula
$$
\Res \biggl[ \sum_i a_i \Lambda^i , \sum_j b_j \Lambda^j \biggr] = \nabla
\sum_k [k](a_kb_{-k}) .
$$

There is a unique linear map
$$
\Res:\Phi_-(\tAA,q)\o_{\tAA}\Om(\tAA)\to\Om(\tAA)/\p\Om(\tAA)
$$
such that $d\Res(A)=\Res(dA)$.

Associated to the equivariant Toda lattice, we have the basic sequence of
functionals
$$
h_n = \frac{1}{n+1} \Res(L^{n+1}) , \quad n\ge0 ,
$$
with differentials $dh_n=\Res(L^ndL)$. In calculating $h_n$, the following
lemma is convenient.
\begin{lemma}
$$
p_0(n+1) = \sum_{k=0}^n [k+1]\bigl( a_{k+1} \, p_k(n) \bigr)
$$
\end{lemma}
\begin{proof}
Applying the operator $\res$ to the equations $L^{n+1}=L\*L^n$ and
$L^{n+1}=L^n\*L$, we see that
\begin{align*}
p_0(n+1) &= \E^{1/2} p_{-1}(n) + \sum_{k=0}^\infty \E^{-k/2} \bigl(
a_{k+1} \, p_k(n) \bigr) , \\
p_0(n+1) &= \E^{-1/2} p_{-1}(n) + \sum_{k=0}^\infty \E^{k/2} \bigl(
a_{k+1} \, p_k(n) \bigr) .
\end{align*}
Taking $\E^{1/2}$ times the second of these equations minus
$\E^{-1/2}$ times the first, we see that
$$
\nabla p_0(n+1) = \nabla \sum_{k=0}^n [k+1]\bigl( a_{k+1} \, p_k(n)
\bigr) ,
$$
and hence, that
$$
p_0(n+1) = \sum_{k=0}^n [k+1]\bigl( a_{k+1} \, p_k(n) \bigr) +
\alpha(p_0(n+1)) .
$$
This proves the lemma, since $\alpha(p_0(n+1))=0$.
\end{proof}
\begin{corollary}
$$
h_n = \sum_{k=0}^n \frac{k+1}{n+1} \tint \bigl( a_{k+1} \, p_k(n) \bigr) \,
dx
$$
\end{corollary}

For example, using the formulas $a_2=q+\t\PP v+z_1$ and
$$
a_3=\t \bigl( \PP \bigl( \tfrac{1}{4} [2]v^2 + q \bigr) - \half v [2]\PP v
\bigr) + \t^2 \PP v - z_1 v + \half z_2 ,
$$
we see that
\begin{align*}
h_0 &= \tint v \, dx , \\
h_1 &= \tint (\half v^2+a_2) dx = \tint (\half v^2+q+\t v+z_1) dx , \\
h_2 &= \tint (\tfrac{1}{3}vp_0(2)+\tfrac{2}{3}(a_2p_1(2))+a_3) dx \\ 
&= \tint (\tfrac{1}{3}v(v^2+[2]a_2)+\tfrac{2}{3}(a_2[2]v)+
\t \bigl( \half v^2 + q - \half v [2]\PP v \bigr) + \t^2 v -
z_1 v + \half z_2) dx \\
&= \tint (\tfrac{1}{3}v^3 + v[2]q + \t \bigl( \half v^2 + q + \half v
[2]\PP v \bigr) + \t^2 v + z_1 v + \half z_2) dx .
\end{align*}

\begin{proposition}
We have $\Res(L^n\,dK)=dH_n$, where
$$
H_n = h_n - \t h_{n-1} + \sum_{k=1}^{n-1} z_k h_{n-k-1} .
$$
\end{proposition}
\begin{proof}
From \eqref{KL}, \eqref{dLs} and \eqref{dlog}, we see that
\begin{align*}
dK &= dL + \t d\ell - \sum_{j=1}^\infty \frac{z_j}{j} dL^{-j} \\
&= dL + \sum_{k=0}^\infty (k+1)^{-1} \ad(L)^k \biggl( \biggl( - \t +
\sum_{j=1}^\infty \tbinom{j+k}{k} z_j L^{-j} \biggr) L^{-k-1}dL
\biggr) .
\end{align*}
Multiplying by $L^n$ and applying $\Res$, all of the terms with $k>0$ drop
out, and we obtain
$$
\Res(L^n\,dK) = \Res \biggl( \biggl( L - \t + \sum_{j=1}^\infty z_j L^{-j}
\biggr) L^{n-1}dL \biggr) ,
$$
which equals $dH_n$.
\end{proof}

Let $\d_v$ and $\d_u$ be the variational derivatives with respect to $v$
and $u=\log(q)$.
\begin{corollary}
We have $\d_vH_n=p_0(n)$, $\d_uH_n=qp_1(n)$, $\d_v\Hbar_n=\pbar_0(n)$ and
$$\d_u\Hbar_n = q\pbar_1(n)-\t\PP\pbar_0(n).$$
\end{corollary}
\begin{proof}
The formulas for $\d_vH_n$ and $\d_uH_n$ follow since
$dK=dv+q\,du\,\Lambda^{-1}$. The formulas for $\d_v\Hbar_n$ and
$\d_u\Hbar_n$ now follow by taking conjugates, bearing in mind that
$\vbar=v-\t\PP u$.
\end{proof}

For example, we have
\begin{align*}
H_0 &= h_0 = \tint v \, dx , \\
H_1 &= h_1-\t h_0 = \tint (\half v^2+q+z_1) dx , \\
H_2 &= h_2-\t h_1+z_1h_0
= \tint (\tfrac{1}{3}v^3 + v[2]q + \half \t v [2]\PP v + 2 z_1 v - \t z_1 +
\half z_2) dx .
\end{align*}

It is now easy to show that the equivariant Toda lattice is Hamiltonian.
Applying $\res$ to the equation $[K,L^n]=\t\p L^n$, we see that
$$
\nabla p_{-1}(n) = \nabla(qp_1(n)) + \t\p p_0(n) .
$$
It follows that $\d_nv=\nabla p_{-1}(n)=\nabla(qp_1(n))+\t\p p_0(n)$. In
conjunction with the formula $\d_nu=\nabla p_0(n)$, we conclude that
$$
\d_n \begin{bmatrix} v \\ u \end{bmatrix} =
\begin{bmatrix} \t\p & \nabla \\ \nabla & 0 \end{bmatrix}
\begin{bmatrix} \delta_vH_n \\ \delta_uH_n \end{bmatrix} .
$$
Since $\dbar_nv=\nabla(q\pbar_1(n))$ and $\dbar_nu=\nabla\pbar_0(n)$, we
also conclude that
$$
\dbar_n \begin{bmatrix} v \\ u \end{bmatrix} = \begin{bmatrix} \t\p &
\nabla \\ \nabla & 0 \end{bmatrix}
\begin{bmatrix} \delta_v \Hbar_n \\ \delta_u\Hbar_n \end{bmatrix} .
$$
In other words, the equivariant Toda lattice is Hamiltonian with respect to
the Hamiltonian structure
\begin{align*}
\{v(x),v(y)\} &= \t \p \delta(x-y) , & \{v(x),u(y)\} &= \nabla_x\delta(x-y)
, & \{u(x),u(y)\} &= 0 .
\end{align*}

The relationship between the equivariant Toda lattice (with $z_k=0$, $k>0$)
and the equivariant Gromov-Witten invariants of $\CP^1$ leads to a new
proof of the Toda conjecture for the (non-equivariant) Gromov-Witten
invariants of $\CP^1$. (See \cite{toda} for a discussion of this conjecture
and further references.) We see that the descendent flow $\p_{k,Q}$ is
the limit of the flow $\frac{1}{(k+1)!}\d_{k+1}$ as $\t\to0$, and hence
has Hamiltonian
$$
\tfrac{1}{(k+1)!} \lim_{\t\to0} h_{k+1} .
$$
Likewise, by \eqref{pP}, the descendent flow $\p_{k,P}$ is the limit of the
flow
$$
\tfrac{1}{(k+1)!} \t^{-1} ( \d_{k+1} - \dbar_{k+1} ) - \tfrac{1}{k!} c_k
(\d_k+\dbar_k)
$$
as $\t\to0$, and hence has Hamiltonian
$$
\lim_{\t\to0} \Bigl( \tfrac{1}{(k+1)!} \t^{-1} ( H_{k+1} - \Hbar_{k+1} ) -
\tfrac{1}{k!} c_k (H_k+\Hbar_k) \Bigr) .
$$
Let $\ell_0$ equal the limit as $\t\to0$ of $\ell$. Since
$L=K-\t\ell=K-\t\ell_0+O(\t^2)$ and
$$
\Lbar = \bar{K}+\t\bar{\ell} = K + \t(\ell_0-\PP u) + O(\t^2) ,
$$
we have
\begin{align*}
\t^{-1} \bigl( H_{k+1} - \Hbar_{k+1} \bigr) &= \t^{-1} \tfrac{1}{k+2}
\Res(L^{k+2}-\Lbar^{k+2}) - \tfrac{1}{k+1} \Res(L^{k+1}-\Lbar^{k+1}) \\
&= \Res(K^{k+1}(\PP u-2\ell_0)) .
\end{align*}
It follows that $\p_{k,P}$ has Hamiltonian $\frac{1}{(k+1)!}
\Res(K^{k+1}(\PP u-2(\ell_0+c_k)))$. An equivalent formula was
conjectured by Eguchi and Yang \cite{EY} and proved by Carlet,
Dubrovin and Zhang \cite{CDZ}.


\begin{thebibliography}{99}

\bibitem{CDZ} G. Carlet, B. Dubrovin and Y. Zhang, to appear.

\bibitem{EY} T. Eguchi and S.-K. Yang, \emph{The topological $\CP^1$ model
and the large-$N$ matrix integral.} Modern Phys. Lett. \textbf{A 9} (1994),
2893--2902. \texttt{<hep-th/9407134>}

\bibitem{toda} E. Getzler, \emph{The Toda conjecture.} In ``Symplectic
geometry and mirror symmetry (KIAS, Seoul, 2000),'' eds. K.~Fukaya et al.,
World Scientific, Singapore, 2001, pp. 51--79. \texttt{<math.AG/0108108>}

\bibitem{eqtoda} E. Getzler, \emph{The equivariant Toda lattice, I.}
\texttt{<math.AG/0207025>}

\bibitem{KZ} B. Khesin and I. Zakharevich, \emph{Poisson-Lie group of
pseudodifferential symbols and fractional KP-KdV hierarchies.}
C. R. Acad. Sci. Paris S\'er. I Math. \textbf{316} (1993),
621--626. \texttt{<hep-th/9311125>}

\bibitem{Kumar} K. Kumar, \emph{Expansion of a function of noncommuting
operators.} J. Math.\ Phys.\ \textbf{6} (1965), 1923--1927.

\bibitem{OP} A. Okounkov and R. Pandharipande, \emph{The equivariant
Gromov-Witten theory of $\mathbf{P}^1$.} \texttt{<math.AG/0207233>}

\bibitem{P} R. Pandharipande, \emph{The Toda equations and the
Gromov-Witten theory of the Riemann sphere.} Lett. Math. Phys. \textbf{53}
(2000), 59--74. \texttt{<math.AG/9912166>}

\bibitem{rahul} R. Pandharipande, private communication (2000).

\bibitem{Zhang} Y. Zhang, \emph{On the $\CP^1$ topological sigma model and
the Toda lattice hierarchy.} J. Geom. Phys. \textbf{40} (2002), 215--232.

\end{thebibliography}
\end{document}